\documentclass{amsart}
\usepackage{amsthm,amsfonts,amsmath,amssymb,latexsym,epsfig,eucal}

\newtheorem{theorem}{Theorem}
\newtheorem{lemma}[theorem]{Lemma}
\newtheorem{corollary}[theorem]{Corollary}

\newenvironment{example}{\medskip \refstepcounter{theorem}
\noindent  {\bf Example \thetheorem}.\rm}{\,}

\def\mb#1{{\mathbb #1}}
\def\mc#1{{\mathcal #1}}

\def\over#1{\overline{#1}}

\def\<{\langle}                                    
\def\>{\rangle}  
\def\ddb{\partial  \overline{\partial}}

\begin{document}

\title{Einstein and scalar flat Riemannian metrics}
\author{Santiago R. Simanca}
\thanks{Supported by the Simons Foundation Visiting Professorship award 
number 657746.}
\address{Department of Mathematics, Courant Institute of Mathematical Sciences,
251 Mercer St., New York, NY 10012} 
\email{srs2@cims.nyu.edu}

\begin{abstract} 
On a given closed connected manifold of dimension two, or greater, we 
consider the squared $L^2$-norm of 
the scalar curvature functional over the space of 
Riemannian metrics of fixed volume. 
We prove that its critical points have constant scalar curvature, and use 
this to show that a metric is a solution of the critical point equation if, 
and only if, it is either Einstein, or scalar flat.
\end{abstract}


\maketitle

\section{Introduction} \label{s1}
The metric $g$ of a closed connected Riemannian manifold $(M^n,g)$ is 
Einstein if it satisfies 
the tensorial relation 
\begin{equation} \label{ei}
r_g =\frac{s_g}{n}g\, , 
\end{equation}
where $r_g$ and $s_g$ are the Ricci and scalar curvature tensors of $g$,
respectively. By using the trace of the differential Bianchi identity, 
this relation implies that $s_g$ must be constant when $n\geq 3$.
If $n=2$, all metrics satisfy this relation. In 
that case, the Einstein condition is strengthened to requiring that $s_g$ 
be constant also. 

It is quite remarkable that there are manifolds that carry 
Einstein metrics with scalar curvatures of opposite signs.
The known examples of such, found by Catanese and LeBrun \cite{cale}, are 
diffeomorphic products of K\"ahler surfaces, and the metrics on the surface 
factors, which are homeomorphic but not diffeomorphic,
are K\"ahler-Einstein relative to complex structures whose corresponding 
first Chern classes are negative, and positive, respectively.    

The extremal metrics of Calabi \cite{ca} are the critical points of the 
squared $L^2$-norm of the scalar curvature functional among 
K\"ahler metrics whose K\"ahler forms represent
a given cohomology class.  If the said class is a multiple of $c_1$ 
(which, therefore, must be a signed class), when $c_1<0$, or when $c_1>0$ and 
the manifold is a complex surface with reductive automorphism group, 
the extremal metric exists, and is Einstein 
\cite{a,yau1}, \cite{tia}. The Einstein metrics of \cite{cale} 
are found appealing to these two results, once  
the complex structures with $c_1<0$ or $c_1>0$ in the factors of the 
underlying manifolds have been identified. (Naturally, these complex 
structures cannot be homotopically equivalent.) But more to the point of 
this note, 
this strongly suggests the general assertion that only Einstein, and scalar flat
 metrics, appear as critical points of the squared $L^2$-norm 
of the scalar curvature functional among Riemannian metrics of fixed volume. 

The critical point equation of this functional, and the 
fact that Einstein and scalar flat Riemannian metrics are solutions of it, have
been known for a while \cite[(5.1), p. 291]{berg}
(see also \cite[Corollary 4.67, p. 133]{be}).
Remarkable geometric analysis based on its use had been
produced up to that point, and ever since, both in the K\"ahler and purely 
real context. So it is slightly surprising that to this date, and except in
some particular cases, it remains unknown if these are the only type of 
critical points that the functional has. In this note, we show that all the 
critical points must have 
constant scalar curvature, from which it follows that they must be Einstein or 
scalar flat metrics, thus filling the void. The characterization 
attained might prove of some use in the study of these metrics. 

If we have a path segment of almost Hermitian structures,
if the metrics at both ends of the path exhibit differences in sign for some 
metric tensor, then somewhere in between some topological condition should 
arise as you cross the point where the sign changes. This simple minded 
thought is our main motivation for attempting to close the gap in the 
description of the critical points of the aforementioned functional.
In light of the result in \cite{cale} alluded to above, 
the Einstein equation is too strong of 
a metric tensor to keep in mind for this idea to work, unless some change is
introduced in such a path by some other method as the intermediate 
metric becomes scalar flat. 

\section{The energy $S(g)$ of a Riemannian metric}
In this section, we prove our results. We make the note self contained 
by deriving relevant equations, after recalling some basic variational formulas.    

We let $M$ be a closed connected manifold of dimension $n$. We denote by
$\mc{M}$ the set of all Riemannian metrics on $M$, and by 
$\mc{M}_1$ the subset of all metrics of
volume one.

In order to treat Riemannian functionals of interest as differentiable 
mappings, the infinite dimensional manifold $\mc{M}_1$ is commonly topologized 
using a Sobolev or $C^{k,\alpha}$ norm of sufficiently high order.
 At $g\in \mc{M}_1$, its tangent space $T_g \mc{M}_1$ consists of the 
space of symmetric $2$-tensors whose traces are $L^2$ orthogonal to the 
constants.  Indeed, if the metric $g$ is deformed 
infinitesimally in the direction of the symmetric two tensor $h$, then the
volume form  $d\mu_g$ varies according to the formula
\begin{equation}
\frac{d}{dt}d\mu_{g+th}\mid_{t=0} = 
\frac{1}{2}{\rm trace}_g\, h \, d\mu_{g} 
\label{vf}
\end{equation}
(see, for instance, \cite[(2.10) p. 288]{berg}), 
and thus, $h$ is in the tangent space if, and 
only if, it satisfies the 
said condition. The $g$-trace of any symmetric two tensor $h$ is just the
pointwise inner-product $(g,h)_g$. We have that
$s_g=(g,r_g)_g$.

We let ${\mathcal S}^{p}(M)$ be the bundle of symmetric 
$p$-tensors on $M$. Covariant differentiation defines a map
$$
\nabla_g^p : {\mathcal S}^{p}(M)
\rightarrow \Omega^{1} M \otimes {\mathcal S}^{p}(M)\, ,
$$ 
which can be composed with the symmetrization operator to define \cite{beeb}
$$
\delta^{*}_g : {\mathcal S}^{p}(M)
\rightarrow {\mathcal S}^{p+1}(M)\, .
$$
The metric dual of $\delta_g^{*}$ defines the operator
\begin{equation}
\delta_g :{\mathcal S}^{p+1}(M) \rightarrow {\mathcal S}^{p}(M)\, .
\label{di}
\end{equation}

The scalar curvature $s_g$ varies according to the expression
\begin{equation}
\frac{d}{dt}s_{g+th}\mid_{t=0}  = 
\Delta_g ({\rm trace}_g\, h) + \delta_g(\delta_g h) -(r_g, h)_g \, ,
\label{vfr}
\end{equation}
where $\Delta_g$ is the Laplacian of the metric $g$
\cite[(2.11) p. 288]{berg}. 
Notice that the trace of the differential Bianchi identity mentioned earlier 
yields the identity $\delta_g r_g=-\frac{1}{2} ds_g$, so if $n>2$ and
(\ref{ei}) holds, the function $s_g$ must be constant. 

The gradient of the Riemannian functional of Hilbert 
\begin{equation}\label{hi}
\begin{array}{rcl}
\mc{M} & \stackrel{H}{\rightarrow} & \mb{R} \vspace{1mm} \\
 g & \rightarrow & {\displaystyle \int_M s_g d\mu_g }\, ,
\end{array}
\end{equation}
is well-known \cite{hilb}. By (\ref{vf}) and (\ref{vfr}),
we obtain that
$$
\frac{d}{dt}H(g+th)\mid_{t=0}=
\frac{d}{dt} \int s_{g+th}d\mu_{g+th}\mid _{t=0} = \int
\left( \frac{s_g}{2}g-r_g, h\right)_{g} d\mu_g\, ,  
$$
and if we restrict the domain of $H$ to $\mc{M}_1$, the ensuing 
Lagrange multiplier allows us to conclude that its critical points are 
the Einstein metrics, those that satisfy the tensorial relation
(\ref{ei}). Notice that this holding for any $g$ when $n=2$ 
is a reflection of the fact that the functional $H$ is then constant,
$4\pi \chi(M)$, by the Gauss-Bonnet theorem. In this case, the 
uniformization theorem 
produces a metric of constant scalar curvature in the conformal class,
an Einstein representative of the class.  

A natural Riemannian functional with more flexible critical points is given by  
\begin{equation}
\begin{array}{rcl}
\mc{M}_1 & \stackrel{S}{\rightarrow} & \mb{R} \\
g & \mapsto & S(g)={\displaystyle \int_M s_{g}^2 d\mu _{g}}\, .
\end{array} \label{func}
\end{equation}
We think of this as the stored energy of $M$ when being in the state defined
by $g$, and ask for the states of $M$ that extremize the energy, 
the critical points of $S$. By the extreme case 
of the Cauchy-Schwarz inequality, if $n>2$, the critical states of 
$H$ are subsumed into the critical states of $S$ 
of constant scalar curvature, which, as we show here, is a property
these latter metrics all have. 

By (\ref{vf}) and (\ref{vfr}), we see that
\begin{equation}
\begin{array}{rcl}
{\displaystyle \frac{d}{dt}\int s^2_{g+th}\mid_{t=0}} & = & {\displaystyle 
\int 2s_g( \Delta_g ({\rm trace}_g\, h) + \delta_g(\delta_g h) -(r_g, h)_g 
+ \frac{1}{4}s_g {\rm trace}_g\, h)d\mu_g } 
\\ & = & {\displaystyle \int \left( \left( 2\Delta_g s_g +
\frac{s_g^2}{2}\right) g +2\nabla_g d s_g -2s_g r_g, h \right)_g d\mu_g }\, .
\end{array} \label{vfr2}
\end{equation}
To our knowledge, this expression was originally obtained by M. Berger
\cite{berg}.

\begin{lemma} \label{un}
\cite[(5.1), p. 291]{berg}, 
\cite[Corollary 4.67, p. 133]{be}.
A Riemannian metric $g \in \mc{M}_1$ is a critical point of
the functional $S$ in {\rm (\ref{func})} if, and only if, we have that
\begin{equation}
\nabla_g S=\left( 2\Delta_g s_g +
\frac{s_g^2}{2}\right) g +2\nabla_g d s_g -2s_g r_g =\lambda_g g\, , 
\label{cm}
\end{equation}
where
\begin{equation}
\lambda_g=\frac{n-4}{2n\mu_g(M)}\int s_g^2 d\mu_g 
=\frac{n-4}{2n }\int s_g^2 d\mu_g \, . \label{co}
\end{equation}
Einstein or scalar flat metrics are critical points.
\end{lemma}

{\it Proof}. By (\ref{vfr2}), $\nabla_g S$ is orthogonal 
to all symmetric two-tensors of trace orthogonal to the constants, so parallel
to $g$. The value of the proportionality constant follows 
by computing the trace of the resulting tensorial identity (\ref{cm}), 
 and integrating the 
resulting functions over $M$ with respect to the measure defined by $g$. 
\qed 
\medskip

Given a Riemannian metric $g$, we denote by $\pi_g$ the $L^2$-projection
operator of $L^2$-functions onto the constants, so in effect,
$S(g)=\pi_g(s_g^2)$.  If $g \in \mc{M}_1$ is
a critical point of the functional $S$, 
the trace of (\ref{cm}) implies the equation
\begin{equation}
(2n-2)\Delta_g s_g +\frac{n-4}{2}\left(s_g^2- \pi_g(s_g^2)\right) = 0
\label{id}
\end{equation}
for the scalar curvature $s_g$ of $g$.
This is precisely the critical point equation of the 
functional $S$ restricted to the set of
volume preserving metrics in the conformal class of $g$. It may be derived
directly from known variational formulas under conformal changes of the
metric \cite[f),g), p. 59]{be}.

As a map defined on the entire set of metrics $\mc{M}$, we have that
\begin{equation} \label{map}
\nabla_g S=\left(\frac{2n-2}{n}\Delta_g s_g +
\frac{n-4}{2n}s_g^2\right)g +2\left(\nabla_g d s_g+\frac{\Delta_g s_g}{n}g-
s_g ( r_g -\frac{s_g}{n}g)\right)\, .  
\end{equation}
Thus, along stationary directions $h \in T_g \mc{M}$ of trace 
orthogonal to the constant, we have that 
$$
\int_M \< \nabla_g S, h\>_g d\mu_g = 0 =  
\frac{n-4}{2n}\pi_g (s_g^2)\int_M (g, h)_g d\mu_g \, , 
$$
and if $g$ is a critical metric of (\ref{func}), by Lemma \ref{un},    
$\nabla_g S= \frac{n-4}{2n}\pi_{g} (s^2_{g})g$. In this latter case,
if $n \neq 4$, and $g$ is not scalar flat, in a sufficiently small 
neihborhood $U_g$ of $g$, the level set 
$\{ \tilde{g}\in U_g: \; S(\tilde{g})= S(g) \} 
\subset \mc{M}$ defines a submanifold of the space of metrics whose tangent 
space at $g$ coincides with $T_g \mc{M}_1$.  

We use these facts to characterize certain directions in the ``light cone'' 
of the Hessian of $S$ at a critical point of (\ref{func}), when these 
directions exist.  
We presume the result to be known, but we have not been able
to locate in the literature a reference to it, or to a more general statement.

\begin{lemma} \label{min}
Suppose that $n\neq 4$, and let $g \in \mc{M}_1$ be a critical metric of the 
functional $S$ in {\rm (\ref{func})}. If $g(t) \in 
\{ \tilde{g}\in U_g: \; S(\tilde{g})= S(g) \} \cap \mc{M}_1$ is a path of
metrics on the $S$ level set of $g$ such that $g(0)=g$, 
then $h=\dot{g}(0)$ is in the kernel of the quadratic form given by the 
second variation of $S$ at $g$ in the direction $h \in T_g \mc{M}_1$, and we 
have that
\begin{equation} \label{seva}
D^2 S_g(h,h)= 0 =  
\frac{n-4}{2n}\pi_g(s_g^2 ) \left( - 
\int (z_g(h),z_g(h))_g d\mu_g +\frac{n-2}{2n}\int (g,h)_g(g,h)_g d\mu_g 
\right)\, , 
\end{equation}
where $z_g(h)$ is the $g$-trace-free component of $h$. 
In particular, if $n\neq 2,4$, and $g$ is a nonscalar flat critical point of 
$S$, if $z_g(h)=0$, or if $z_g(h)=h$, $h$ must be identically zero.
\end{lemma}

{\it Proof}. Since $S(g(t))=S(g)$, we have that
$$
S(g(t))-S(g)=0=\int_0^t \nabla_{g(s)}S \cdot \dot{g}(s) ds =\int_0^t
\left( \int_{M} (\nabla_{g(s)}S, \dot{g}(s))_{g(s)} d\mu_{g(s)}\right)ds  \, ,
$$
and we can find a sequence $t_n \searrow 0$, such that
$$
\int_{M} (\nabla_{g(t_n)}S, \dot{g}(t_n))_{g(t_n)} d\mu_{g(t_n)} =0 \, .
$$
Since $\dot{g}(t)$ is of $g(t)$-trace that is orthogonal to the constants, 
we have that along the sequence of stationary directions 
$\dot{g}(t_n)$ at $g(t_n)$ for the functional $S$, we have that
$$
(D_{g(t_n)}S)(\dot{g}(t_n))=0=\frac{n-4}{2n}\pi_{g(t_n)} (s^2_{g(t_n)})
\int (g(t_n), \dot{g}(t_n))d\mu_{g(t_n)}\, ,
$$
and $\nabla_{g(0)}S=\frac{n-4}{2n}\pi_{g(0)} (s^2_{g(0)})
g(0)$. Thus, the quadratic form given by the Hessian 
$D^2_{g}S(k,h)=\frac{d}{dt}(D_{g(t)}S)(h)\mid_{t=0}$ at $k=h$, 
which vanishes, can be computed by computing the directional derivative of
$\frac{n-4}{2n}\pi_{g} (s^2_{g}) \< g, h\>_g$ along the path $g(t)$ evaluated
at $t=0$. 

We write $( \nabla_g S, h )_g = 
(\nabla_g S)_{ij}h_{pq}g^{ip}g^{jq}=
((n-4)/2n)\pi(s_g^2)g_{ij}h_{pq}g^{ip}g^{jq}$. 
Then, computing the Hessian quadratic form at $h$ by moving along the path
$g(t)$, by the chain rule we see that  
$$ 
\begin{array}{rcl}
0=D^2 S_g(h,h) & =  & {\displaystyle \frac{n-4}{2n} \pi_g(s_g^2) \left(        
-\int (h,h)_g d\mu_g +\frac{1}{2}\int (g,h)_g(g,h)_g d\mu_g -
\left(\frac{2}{n}\right)\int (g,h)_g d\mu_g 
\int (g,h)_g d\mu_g \right) } \vspace{1mm} \\ 
& = & {\displaystyle \frac{n-4}{2n} \pi_g(s_g^2) \left(        
-\int (h,h)_g d\mu_g +\frac{1}{2}\int (g,h)_g(g,h)_g d\mu_g 
\right)} \, ,  
\end{array}
$$ 
where the first and last term on the right side of the first equality 
are the result of differentiating the three $g$s and $\pi(s_g^2)$ terms in 
$\frac{n-4}{2n}\pi_{g} (s^2_{g})(g,h)$, respectively, the second equality
being then clear since $h$ has trace orthogonal to the constants. 

If we write symmetric tensors $h,k$ in terms of their trace and trace-free 
components, we obtain that
$$
(h,k)_g=\left( \frac{(h,g)_g}{n}g+ z_g(h), 
 \frac{(k,g)_g}{n}g+ z_g(k)\right)_g =
 \frac{1}{n}(h,g)_g (k,g)_g  + (z_g(h),z_g(k))_g \, . 
$$ 
We use this fact in the preceeding expression with $k=h$. 
The asserted form for $D^2 S_g(h,h)$ follows. 
 \qed
\medskip

There is ample amount of remarkable research involving the functional $S$, but  
its complete citation here is out of the scope of the note. 
Insofar as our concern goes, we find that the critical points of $S$ 
have been described in particular cases only (see \cite[page 133]{be},
\cite[$n$=2]{ca}, \cite[Proposition 1.1]{an} or 
its more recent extension in \cite{cat} in the nonnegative case). We proceed to 
show that these critical metrics all lie within the set of metrics of constant 
scalar curvature, from which point on, their final description as Einstein  
or scalar flat metrics follows.

It is not immediate that the scalar curvature $s_g$ of a metric 
satisfying the critical equation (\ref{cm}) must be a constant, but there are
precedents when the dimension of $M$ is either $2$ or $4$. In the latter case, 
$s_g$ is proved to be constant \cite{be} as an application of the 
maximum principle, and as such, the argument has a somewhat ``local nature'' 
flavor (in this regard, Besse's main interest seems to have been on 
manifolds of dimension $4$, perhaps the reason why his result was only 
derived in that case). The argument used in dimension $2$ has, by 
contrast, a ``global'' flavor built into it 
(as does the nonnegative case 
in \cite[Proposition 1.1]{an}, Anderson's interest 
on three manifolds surely the reason why the proof is written
for $n=3$ only). We pause to analyze the details of these two arguments.

When $n=4$, a solution $g$ to (\ref{id}) must be such that $s_g$ is a 
harmonic function, and therefore, constant by the maximum principle
\cite[Corollary 4.67]{be}. 

In general, a solution of (\ref{id}) is constant
if, and only if,
$s_g^2 -\pi_g(s_g^2)$ has a zero of order at least three, 
in which case $s_g =\pi_g (s_g)$ equals one of the square roots of
$\pi_g(s_g^2)$. 
 For (\ref{id}) implies that this zero must be a 
zero of $s_g \mp \sqrt{\pi_g(s^2_g)}$ of infinite order, and the vanishing of 
this function would then be a consequence of Aronszajn's unique continuation 
theorem for solutions to elliptic equations of order two \cite{ar}. 
Alternatively, 
if there exists a point in $M$ where the function $u_g=s_g-\pi_g(s_g)$ 
vanishes to order at least three, then $s_g$ is constant, and $u_g \equiv 0$
so the said point is in fact a zero of infinite order. Indeed, in terms of 
the function $u_g$, (\ref{id}) is given by the equivalent expression
$$
(2n-2)\Delta_g u_g +\frac{n-4}{2}(u_g^2+2u_g \pi_g(s_g)+(\pi_g(s_g))^2-\pi_g
(s_g^2)) = 0 \, ,
$$
which evaluated at the point implies that the constants $\pi_g (s^2_g)$ 
and $(\pi_g(s_g))^2$ are equal. The assertions then follow by
the extreme case of the Cauchy-Schwarz inequality.

But if $s_g$ is a solution to (\ref{id}), the existence of a point where 
$s_g^2 -\pi_g(s_g^2)$ has a zero of order three is hard to prove on its own, 
if at all possible. Although the function $s_g$ is constant, as we shall see,
the proof of this fact requires information contained in equation (\ref{cm})
itself, which is somewhat lost in passing to its trace equation (\ref{id}).
This subtle point is illustrated well in the known proof that 
$s_g$ is constant when $n=2$, an argument that,  
unexpectedly, is more elaborate than the $n=4$ case above. 

Indeed, let $M$ be a differentiable surface. By passing to a double 
covering if necessary, let us assume that $M$ is oriented. Then $M$ can be 
provided with a compatible complex structure $J$, which is 
defined by taking an orthonormal $g$-frame $\{ e_1, e_2\}$ and declaring that
$Je_1:=e_2$. This makes of $M$ a complex manifold of dimension $1$, and 
the metric $g$ is K\"ahler 
with K\"ahler form $\omega_g(\, \cdot \, ,\, \cdot \,)= 
g(J\, \cdot \, ,\, \cdot \,)$. 
Since $g$ is a critical metric of the functional $S$, it must be 
critical also for the restriction of this functional to the
submanifold of K\"ahler metrics that represent the 
cohomology class $[\omega_g]$, or extremal in the sense of Calabi 
\cite{ca}. (In this dimension, this is the same as metrics in the conformal 
class of $g$.) Thus, if  $\partial_g^{\#}f$ is the vector field 
 defined by the identity
$g(\partial^{\#}f, \, \cdot \,)=\over{\partial} f$ \cite{lesi3}, 
we must have that $\partial_g^{\#}f$ is holomorphic \cite{ca}.

For dimensional reasons, we have that
$\rho_g=\frac{s_g }{2}\omega_g$, where $\rho_g$ is the Ricci form of $g$,  
and the identity $-2i\ddb f =(\Delta_g f) \omega_g$ holds for any function $f$. 
By  (\ref{id}), we obtain that 
$\Delta_g^2 s_g=s_g\Delta_g s_g- (\nabla^g s_g, \nabla^g s_g)$, and so
$$
\Delta_g^2 s_g + 4(\rho_g, i\ddb s_g)=-(\nabla^g s_g, \nabla^g s_g)\, ,
$$
an equation that we rewrite as
$$ 
4(\over{\partial} \partial^{\#})^{*}\over{\partial} \partial^{\#}s_g =
\Delta_g^2 s_g + 4(\rho, i\ddb s_g)+(\nabla^g s_g,\nabla^g s_g)=0\, .
$$ 
Thus, $\partial_g^{\#}s_g$ is holomorphic. 
If $M$ is either a hyperbolic or a parabolic Riemann surface,  
$\partial_g^{\#}s_g =0$, and $s_g$ is constant.  
In the elliptic case, $\partial_g^{\#}s_g=0$ also, but for a different 
and less elementary reason: The Kazdan-Warner invariant \cite{kw} vanishes, 
and this measures the obstruction of $g$ to being 
conformally equivalent to the standard metric. It follows that $g$ itself is 
the standard metric, and $s_g$ is constant. Thus, for any closed 
connected surface $M$, 
a critical point of (\ref{func}) must be a metric of constant scalar curvature
(as proved in \cite{ca} in the oriented case; the positive case of this 
proof covers the like claimed but unproven statement for $n=2$ in \cite{cat}).

The point of this argument 
is that by the dimensional identity $r_g=\frac{s_g}{2}g$, if $g$ is a 
critical point of the functional $S$, the $g$-trace-free component of
$\nabla^g d s_g$ vanishes, and (\ref{cm}) and (\ref{id}) are then equivalent 
equations implying that $\nabla^g s_g$ is a conformal vector
field (that if $M$ is oriented, and $g$ is a K\"ahler metric relative
to the compatible complex structure, is the same as saying that  
the vector field $\partial^\# s_g$ is holomorphic). Now for any conformal 
vector field $X$ of a Riemannian metric $g$, the scalar curvature
$s_g$ satisfies the Kazdan-Warner identity 
$$
\int X (s_g) d\mu_g = 0\, ,  
$$
\cite[Theorem II.9]{boez}, and
this applied with $X=\nabla^g s_g$ to our critical metric $g$ implies that 
$s_g$ is constant.  
However, it is not the case that all solutions of 
(\ref{id}) (as an equation in $s_g$) are constants. For the ordinary 
differential equation 
$$
2 \ddot{u} = c^2- u^2, \quad \text{$c$ constant,}
$$
has periodic solutions, and, for instance, on flat tori of all 
dimensions, the partial differential equation 
$$
2 \Delta_g u -(u^2- c^2)=0
$$
has many nonconstant solutions, in addition to $u=c$. On a surface,
by using the Kazdan-Warner identity above, we see that the critical point 
equation of $S$ singles out 
the metric whose scalar curvature is the constant solution to (\ref{id})
in the set of all metrics of area one in the conformal class, the 
uniformization theorem for the conformal class that the metric defines.

This last point indicates that the maximum principle alone in (\ref{id}) cannot
be enough to prove that a critical metric of $S$ has constant scalar curvature.
However, we have the following partial result:

\begin{lemma}\label{le2}
Consider a critical metric $g$ of the functional
$S$ in {\rm (\ref{func})} with scalar curvature $s_g$. Then:
\begin{enumerate}
\item For $n> 4$, if $\max{s_g} \geq |{\min{s_g}}|$, then $s_g$ must be
a nonnegative constant.
\item For $n < 4$, if $\min{s_g} \leq -|{\max{s_g}}|$, then $s_g$ must be
a nonpositive constant.
\end{enumerate}
\end{lemma}

{\it Proof}.  (1) Since $\max{s_g} \geq |{\min{s_g}}|$, a
maximum of $s_g$ is also a maximum of $s_g^2$.
At a maximum $p$ of $s_g$, $\Delta_g s_g \geq 0$. If
$s^2_g(p)  > \pi_g(s_g^2)$, the value of the left side of (\ref{id}) at $p$ 
would be positive. Thus, 
the maximum of $s^2_g$ must equal its average $\pi_g(s_g^2)$, so 
the function $s^2_g$ is constant, and $s_g$ itself is constant,
nonnegative.  
 
(2) At a minimum $p$ of $s_g$, $\Delta_g s_g \leq 0$. If
$\min{s_g} \leq  -|{\max{s_g}}|$, then $p$ is a maximum of
the function $s_g^2$. If 
$s^2_g(p)  > \pi_g(s_g^2)$, the value of the left side of (\ref{id}) at $p$ 
would be negative. Thus,  
the maximum of $s^2_g$ equals its average  
$\pi_g(s_g^2)$, so $s_g^2$ is constant, and the function $s_g$ itself 
is constant, nonpositive.
\qed

Notice that below the threshold $n=4$, this result handles the
generic case, and leaves unsettled the nongeneric one, exactly the opposite
of what happens above it.
 
Our proof that the critical points of $S$ have constant scalar 
curvature involves a global consideration, as is the case of the $n=2$ 
argument presented earlier. This global consideration is of 
a slightly different nature, and applies to resolve the situations complementary
to those covered by the maximum principle in Lemma \ref{le2}. We shall make use
of the content of Lemma \ref{min}, with the restrictions it imposes on 
stationary directions at a critical metric.
Notice that along any symmetric tensor $h \in T_g \mc{M}_1$ for which 
$z_g(h)=0$, and any constant $c$, 
by (\ref{map}), we have that 
\begin{equation} \label{fact}
\< \nabla_g S, h\>= 
\int \left(\frac{2n-2}{n}\Delta_g s_g +
\frac{n-4}{2n}(s_g^2 -c^2)\right)(g,h)_g d\mu_g \, . 
\end{equation}
We exploit this fact with $c=\pi_g (s_g)$ below.

\begin{theorem}\label{tlk}
The scalar curvature $s_g$ of any critical metric $g$ of the functional
$S$ in {\rm (\ref{func})} is constant.
\end{theorem}

{\it Proof}. We assume firstly that $n>4$. 
If $\max{s_g} \geq | \min{s_g}|$, the desired conclusion follows by 
Lemma \ref{le2}. So let us assume that $\max{s_g} < | \min{s_g}|$. 
Then the range of $s_g$ must contain negative values, and we must have that 
$$
\min{s_g} \leq -|\max{s_g}| \, .
$$

If equality holds, by the maximum principle on (\ref{id}), 
we conclude that $s_{g}=-\sqrt{\pi(s^2_g)}=\pi(s_g)$. We assume the strict 
inequality, and so $s_g$ is not constant.
 In this case, a maximum of $s_g$ is not necessarily
related to an extremal of $s_g^2$, but at such a point $p$, by the 
maximum principle on (\ref{id}), we conclude that $s_g^2(p)\leq \pi(s_g^2)$. 
On the other hand, a minimum $p$ of $s_g$ is a maximum of $s_g^2$, 
and by the maximum principle on (\ref{id}), we must have that 
$s_g^2(p) > \pi(s_g^2)$. We rule this situation out as follows.

We consider a perturbation of the Yamabe flow of Hamilton \cite{haY}
given by 
\begin{equation} \label{flow} 
\frac{dg}{dt}=( \pi_g(s_g)-s_g -v )g \, , 
\end{equation} 
where the perturbation term $v$ is a $t$-dependent 
 elliptic pseudodifferential operator 
of order $-2$ on $g$ that is orthogonal to the constants, and is such that
$v(0)=0$. We describe its defining property next, but observe right away
that the solution $g(t)$ of this flow equation that starts at the 
critical metric $g$, while it exists,
produces a path of volume one metrics on the level set $S(g(t))=S(g)$ of $g$.

For notational convenience, we set $u= \pi_g(s_g) -s_g$, $h=u-v$, and also, 
$a_n=2n-2$ and $b_n = \frac{1}{2}(n-4)$, respectively.
 We apply (\ref{fact}) with 
$c=\pi_g(s_g)$. Then we have that  
$$
\begin{array}{rcl} 
\< \nabla_g S, h g\> 
& = & {\displaystyle \int (a_n\Delta_g s_g + b_n(s^2_g- (\pi_g(s_g))^2))h 
d\mu_g }\\ 
& = & {\displaystyle \int (s_g + \pi_g(s_g))( a_n\Delta_g h + b_n(s_g-\pi_g
(s_g))h) d\mu_g }  \, ,   
\end{array}
$$
and so for this to be zero, the function 
$$
(a_n\Delta_g h + b_n(s_g-\pi_g (s_g))h=  
a_n\Delta_g u + b_n(s_g-\pi_g (s_g))u 
- (a_n\Delta_g v + b_n(s_g-\pi_g (s_g))v) 
$$
must be orthogonal to $s_g + \pi_g (s_g)$ for all $t$. 

Since the starting metric of the flow
is a critical point of (\ref{func}), at $t=0$ we have that  
$a_n\Delta_g u + b_n(s_g-\pi_g (s_g))u$ is $L^2$-orthogonal to 
$s_g +\pi_g(s_g)$, and so the desired condition holds then if $v(0)=0$.
We ensure that the condition holds for all time by defining $v$ to 
cancel out the 
projection of $a_n\Delta_g u + b_n(s_g-\pi_g (s_g))u$ onto $s_g +\pi_g(s_g)$ 
as $t$ varies. Thus, we choose $v$ to be the solution of the equation 
\begin{equation} \label{eq14} 
\left( a_n\Delta_g  + b_n(s_g-\pi_g (s_g)\right) v= \Lambda (s_g + \pi_g (s_g))
\, ,   
\end{equation}
where $\Lambda$ is the $t$-dependent constant
$$
\Lambda =\frac{\<s_g +\pi_g(s_g),a_n\Delta_g  + b_n(s_g-\pi_g (s_g))u\>}
{ \| s_g + \pi_g(s_g)\|^2} \, .
$$
Since the initial metric $g$ has nonconstant scalar curvature $s_g$, 
and $b_n \neq 0$, at least for sufficiently small $t$, (\ref{eq14}) can be 
solved for $v$ with $v$ orthogonal to the constants. We have that
$\Lambda(0)=0$, and so this $v(t)$ vanishes when
$t=0$.  
 
Since the regular Yamabe flow admits global solutions in time 
(a result of R. Hamilton that was never published), 
we may use a perturbation argument to show that the initial value problem 
for (\ref{flow}) starting at $g$ has a solution $g(t)$, which is defined
for sufficiently small values of $t$. But we have that
$\dot{g}(0)=(\pi_g (s_g) - s_g)g$ is a $g$-trace free symmetric tensor, 
so by Lemma \ref{min}, we conclude that
$$
D^2 S_g(h,h)= 0=  
\frac{(n-4)(n-2)}{(2n)^2}\pi_g(s_g^2 ) 
\int ( \pi_g (s_g) -s_g)^2 d\mu_g \, .  
$$
This contradicts the fact that $s_g$ is not the constant function.   

If $n=3<4$, by Lemma \ref{le2},
we reduce the analysis to the case where
$$
\max{s_g} \geq |\min{s_g}| \, .
$$
If $s_g \geq 0$, the desired conclusion was proven by Anderson
\cite[Proposition 1.1]{an}. For the remaining possibilities, 
$s_g$ is nonconstant and changes sign, a situation that we rule out
proceeding verbatim as we did in the case above. For the properties of 
the perturbed Yamabe flow that we needed then to show the contradiction were
only conditional on the dimension $n$ being any other than $2$ or $4$. Thus, 
applying Lemma \ref{min} to 
the solution of the perturbed Yamabe flow equation now for $n=3$, we see
that the light cone direction $(\pi_g (s_g) - s_g)g$ of 
the Hessian of $S$ at $g$ has trace of $L^2$-norm zero, once again 
contradicting the fact that $s_g$ is not constant.  
\qed

\begin{corollary} \label{th3}
If $n>2$, the Euler-Lagrange equation of $S$  
is equivalent to the equation 
\begin{equation} \label{te}
2s_g\left( \frac{s_g}{n} g-r_g \right) = 0 \, .
\end{equation}
A metric $g$ is critical if, and only if, 
it is either Einstein, or scalar flat.
\end{corollary}

{\it Proof}. A solution $g$ of the critical point
equation (\ref{cm}) has constant scalar curvature $s_g$, and so 
(\ref{cm}) simplifies to the tensorial relation (\ref{te}). On the other 
hand, if $n>2$, any solution of (\ref{te})  
has constant scalar curvature $s_g$, and so $g$ solves (\ref{cm}) also. 

Let $g$ be a critical metric, so $s_g$ is constant. If $n=2$, 
we always have that $r_g=\frac{s_g}{2}g$, and $s_g$ being constant makes of
$g$ an Einstein metric. 
If $n>2$, and the constant $s_g$ is nonzero, by  
(\ref{te}) we have that $r_g=\frac{s_g}{n}g$, and $g$ is Einstein.
On the other hand, any scalar flat metric, Ricci flat or otherwise, 
satisfies (\ref{te}).   
\qed

If in passing from the Hilbert functional $H$ to the functional $S$ the set 
of Einstein critical points is unaltered, it is conceivable that we 
could go the other way around, and recover at least certain special Einstein 
critical metrics for $S$ as critical points of $H$ of some sort. 
In the K\"ahler context, this program can be carried out for any 
Calabi extremal Einstein 
metric \cite{ca} by using the modified Hilbert functional of
Perelman \cite{pe} over a suitable domain of definition \cite{si2}. In fact, 
it is  possible to see any extremal K\"ahler metric as a critical point of 
this modified functional over a convenient domain, and realize the extremal 
flow \cite{si3} as the flow of the metric along its gradient.

But of more importance, equation (\ref{cm}) may be used 
to study the plausible existence of Einstein metrics 
on a given manifold by a procedure that parallels that for
the plausible existence of strongly extremal K\"ahler 
metrics \cite{si4} (see \cite{sist,sist2,sito} also). We could
first find conformal classes of metrics $[g]$ for which there exists a 
representative $g$ such that $z_g(\nabla^g ds_g
-s_gr_g)=0$, and then within those classes, study the metrics that have
constant scalar curvature. The latter set would contain an Einstein metric,
if any exists. By employing the techniques of Gao and Yau \cite{gaya}, and 
Lohkamp \cite{lohk}, it seems possible to treat the first of 
these problems in the generic case, and perhaps find suitable candidates 
of conformal classes that are represented by metrics of 
negative Ricci curvature. The Einstein metrics, if any, 
 would then be within these classes.

\section{Some remarks on the critical values of $S$}
In dimension two, the set of critical values of $S$ can be 
described completely, a direct consequence of the  
uniformization and Gauss-Bonnet theorems.

\begin{theorem} \label{twod}
For any closed connected surface $M$, the functional $S$ 
in {\rm (\ref{func})} has
exactly one critical value given by $16\pi^2 \chi(M)^2$, $\chi(M)$ the Euler 
characteristic of $M$. In each conformal class of metrics on $M$, 
there exists an Einstein metric $g$ of scalar curvature 
$s_{g}= 4\pi \chi(M)$ that achieves it, unique up to isometries. 
\end{theorem}

A number of results indicate that this theorem has no strong counterpart in 
higher dimension. The entire picture is far from being clear, with the 
nongeneric positive case already 
being substantially more complicated, the generic 
negative case expected to be more so. For emphasis, we summarize some of these 
known facts in the form of a lemma, and briefly discuss their proofs,
for completeness. The statements can be sharpened for $n=3$, after 
the resolution of the geometrization conjecture. We opt for not doing so, as
that would take us immediately far from the scope of the note.
    
\begin{lemma} \label{p5}
If $M$ carries critical metrics  of the functional $S$ in 
{\rm (\ref{func})} of nonnegative 
scalar curvature, then the set of scalar curvatures of all such is bounded 
above. If $n\geq 3$, this set could contain infinitely many elements.   
If $M$ carries a metric $g$ of nontrivial nonnegative scalar curvature, 
and $n\geq 3$, then $0$ must be a critical value of $S$.
\end{lemma}

{\it Proof}. If ${\displaystyle r_g =\frac{s_g}{n}g \geq  0}$, then by Bishop 
comparison theorem \cite{bi}, we have that 
$$
\mu_g (M) \leq \omega_n ({\rm diam}(M,g))^n\, , 
$$
where $\omega_n = \frac{\pi^{\frac{n}{2}}}{\Gamma\left( \frac{n}{2}+1
\right)}$,
$\Gamma$ the Gamma 
function. 
If $s_g > 0$, by Myers' theorem \cite{my}, we have that 
$$
{\rm diam}(M,g) \leq \frac{\pi}{\sqrt{ \frac{s_g}{n(n-1)}}} \, , 
$$
and so $s_g$ must be bounded above since $\mu_g(M)=1$.

On  $M=\mb{S}^2 \times \mb{S}^3$, Wang and Ziller \cite{wazi} construct 
a countably infinite set of volume one Einstein metrics
$g_n$ whose scalar curvatures $s_{g_n}$ are positive, and such that
$s_{g_n}\searrow 0$. In the limit, the metrics collapse.  

The latter assertion is \cite[Theorem 4.32 (ii)]{be}. Since the
conformal Laplacian $L_g= 4\frac{n-1}{n-2}\Delta_g + s_g$
is a strictly positive self-adjoint operator, 
the Yamabe conformal invariant 
$\mc{Y}(M,[g])$ is strictly positive. On the other hand, $M$ carries a metric 
$\tilde{g}$ with negative scalar curvature (see, for instance, \cite{av} and 
\cite{au}), and so $\mc{Y}(M,[\tilde{g}])<0$.
By continuity of the Yamabe invariant, along the segment $(1-t)g+t\tilde{g}$
there exists a metric whose Yamabe invariant is zero. Its conformal class 
contains a metric of zero scalar curvature. This metric is a critical point of 
$S$ of critical value $0$. 
\qed

There is no reason to believe that the scalar flat metric on $M$ 
in the very last argument above is Ricci flat. In fact, if $\Sigma$ is a
complex curve of genus at least two, the blow-up of
$\mb{P}^1(\mb{C})\times \Sigma$ at sufficiently many points 
carries scalar-flat K\"ahler metrics that are never Ricci flat, and
this implies that it carries also K\"ahler metrics of constant scalar 
curvature $c$ for arbitrary $c\in \mb{R}$ 
\cite[Proposition 5, Corollary 1]{lesi1}. Notice that any
product $(M_1, g_1)\times (M_2,g_2)$ of constant scalar curvature 
metrics $g_1$, and $g_2$, such that $s_{g_1}=-s_{g_2}\neq 0$, is  
scalar flat, but not Einstein.

We close the note by discussing a path of Hermitian deformations of  
an Einstein metric on a closed manifold along which, the metrics in the path
have different signs for its
scalar and $J$ scalar curvatures \cite{risi}, a bit intended to justify 
the comment we made at the end of \S \ref{s1}.

\begin{example}
The Calabi-Eckmann manifold $\mb{S}^{2n+1}\times \mb{S}^{2n+1}$ is the total
space of a flat torus bundle over $\mb{P}^n(\mb{C}) \times \mb{P}^n(\mb{C})$, 
and thus, carries an integrable almost complex structure $J$ 
compatible with the standard product metric $g_n$ \cite{caec}. Suppose that 
this metric is dilated in the vertical directions by a factor 
$\varepsilon^2$. Then, the
scalar and $J$-scalar curvature of the dilated metric $g_n^{\varepsilon}$
are given by 
$$
\begin{array}{rcl}
s_{g_{n}^{\varepsilon}}
 & = & 4n(2n+1)+ 4n(1-\varepsilon^2)\, ,\vspace{1mm} \\
s^{J}_{g_{n}^{\varepsilon}} & = & 4n+ 
4n(2n+1)(1-\varepsilon^2)\, , 
\end{array}
$$ 
respectively \cite{dusi}. Notice that $g_n=g_n^1$.

It follows that as $\varepsilon \searrow 0$, the Einstein metric $g_n$ is
deformed and collapses in the Gromov-Hausdorff sense to the product of the
standard Fubini-Study metrics on $\mb{P}^n(\mb{C}) \times \mb{P}^n(\mb{C})$. 
On the other hand, if we blow up the torus fiber letting
$\varepsilon \nearrow \infty$, the Einstein metric $g_n$ transitions 
smoothly from metrics such that 
$s_{g_{n}^{\varepsilon}}>0$ and $s^J_{g_{n}^{\varepsilon}}>0$ 
to ones where both of these scalar tensors have negative values, the
inequality $s_{g_{n}^{\varepsilon}}>s^J_{g_{n}^{\varepsilon}}$ holding always.
\end{example}


\begin{thebibliography}{XX}

\bibitem{an}
M. Anderson, {\it Extrema of curvature functionals on the space of metrics
on $3$-manifolds}, Calc. Var. and P.D.E. 5 (1997), pp. 199-269.
\bibitem{ar}
N. Aronszajn, {\it A unique continuation theorem for solutions of elliptic
 partial differential equations of second order}, J. Math. Pures Appl.,
36 (1957), pp. 235-249.
\bibitem{a} 
T. Aubin, {\it Equations du Type Monge-Amp\`{e}re sur les Vari\'{e}t\'{e}s 
K\"ahl\'eriennes  Compactes},  C. R. Acad. Sci. Paris 283A (1976), pp. 119-121. 
\bibitem{au}
T. Aubin, {\it M\'etriques riemannienen et courbure}, J. Diff. Geom. 11 (1976),
pp. 573-598.
\bibitem{av}
A. Avez, {\it Valeur moyennee du scalaire de courbure sur une vari\'et\'e
compacte}, Applications relativistes, C.R. Acad. Sci. Paris, 256 (1963), pp.
5271-5273.
\bibitem{berg}
M. Berger, {\it Quelques formules de variation pour une structure 
Riemannienne}, Ann. scient. \'Ec. Norm. Sup. Paris, 3 (1970), pp. 285-294.
\bibitem{beeb}
M. Berger \& D. Ebin, {\it Some decompositions of the space of symmetric tensors
on a Riemannian manifold}. J. Diff. Geom. 3 (1969), pp. 379-392.
\bibitem{be}
A. L. Besse, Einstein manifolds, Ergebnisse der Mathematik und ihrer
Grenzgebiete; 3 Folge, Band 10, Springer-Verlag, 1987.
\bibitem{bi}
R. Bishop, {\it A relation between volume, mean curvature and diameter}, 
Not. Amer. Math. Soc, 10 (1963), p. 364.
\bibitem{boez}
J.P. Bourguignon \& J.P. Ezin, {\it Scalar curvature functionals in a conformal
class of metrics and conformal transformations}, Trans. A.M.S. 301 (1987), 
pp. 723-736.  
\bibitem{ca}
E. Calabi, {\it Extremal K\"ahler metrics}, in Seminar on Differential
Geometry (S. T. Yau Ed.), Annals of Mathematics Studies, Princeton University
Press, 1982, pp. 259-290.
\bibitem{caec}
E. Calabi \& B. Eckmann, {\it A class of compact, complex manifolds which
are not algebraic}, Ann. of Math., 58 (1953), pp. 494-500.
\bibitem{cale}
F. Catanese \& C. LeBrun, {\it On the scalar curvature of Einstein manifolds}.
 Math. Res. Lett. 4 (1997), pp. 843–854. 
\bibitem{cat}
G. Catino, {\it Critical metrics on the $L^2$-norm of the scalar curvature},
Proc. Amer. Math. Soc. 142 (2014), pp. 3981–3986.
\bibitem{risi}
H. del Rio \& S.R. Simanca, {\it The Yamabe problem for almost Hermitian 
manifolds}.  J. Geom. Anal. 13 (2003), no. 1, pp. 185-203.
\bibitem{dusi}
C. Dur\'an \& S.R. Simanca, {\it Hermitian metrics on Calabi-Eckmann manifolds}.
Differential Geom. Appl. 17 (2002), no. 1, pp. 55-67. 
\bibitem{gaya}
L.Z. Gao \& S.T. Yau, {\it The existence of negatively Ricci curved metrics on
three manifolds}. Inv. Math. 85 (1986), pp. 637-652.
\bibitem{haY}
R. Hamilton, {\it The Ricci flow on surfaces}, Mathematics and General 
Relativity,  Contemporary Math. 71 Amer. Math. Soc., Providence, RI, 1988, 
pp. 237-262.
\bibitem{hilb}
D. Hilbert, {\it Die Grundiagen der Physik}.  Nachr. Ges. Wiss. G\"ott., 
(1915), p.  395-407.
\bibitem{kw}
J.L. Kazdan \& F.W. Warner, {\it Curvature functions for compact two
manifolds}, Ann. of Math., 99 (1974), pp. 14-47.
\bibitem{lesi1}
C. LeBrun \& S.R. Simanca, {\it On K\"{a}hler Surfaces of Constant Positive
Scalar Curvature}, J. Geom. Anal. 5 (1995), no. 1, pp. 115-127.
\bibitem{lesi3}
C. LeBrun \& S.R. Simanca, {\it On the K\"{a}hler Classes of Extremal
Metrics}, Geometry and Global Analysis (Sendai, Japan 1993),
First Math. Soc. Japan Intern. Res. Inst. Eds. Kotake, Nishikawa \& Schoen.
\bibitem{lohk}
J. Lohkamp, {\it Metric of negative Ricci curvature}, Ann. of Math. 140 
(1994), pp. 655-683.
\bibitem{my}
S.B. Myers, {\it Riemannian manifolds with positive mean curvature}, 
Duke Math. J. 8 (1041), pp. 401-404.
\bibitem{pe}
G. Perelman, {\it The entropy formula for the Ricci flow and its geometric
applications}, preprint DG/0211159, 2002.
\bibitem{si3}
S.R. Simanca, {\it Heat Flows for Extremal K\"ahler Metrics}.
Ann. Sc. Norm. Sup. Pisa Cl. Sci., 5 (4) (2005), no. 2, pp 187-217.
\bibitem{si2}
S.R. Simanca, {\it Hilbert-Perelman's functional and Lagrange multipliers}, 
Proc. Amer. Math.  Soc. 140 (2012), no. 12, pp. 4309–4318.
\bibitem{si4}
S.R. Simanca, {\it Strongly Extremal K\"ahler Metrics}, Ann. Global Anal. Geom.
18 (2000), no. 1, pp. 29-46. 
\bibitem{sist}
S.R. Simanca \& L.D. Stelling, {\it Canonical K\"ahler classes}. Asian J. 
Math.  5 (2001), no. 4, pp. 585-598.
\bibitem{sist2}
S.R. Simanca \& L.D. Stelling, {\it The Dynamics of the Energy of a 
K\"ahler class}, Comm. Math. Phys. 255 (2005), pp. 363-389.
\bibitem{sito}
S.R. Simanca \& C.W. Tonnesen Friedman, {\it The energy of a K\"ahler class 
on admissible manifolds}, Math. Annalen, 31 (2011), pp. 805-834. 
\bibitem{tia} 
G. Tian, {\it On Calabi's Conjecture for Complex Surfaces with Positive 
First Chern Class}, Inv. Math., 101 (1990), pp. 101-172.
\bibitem{wazi}
M. Wang \& W. Ziller, {\it Einstein Metrics on Principal Torus Bundles},
J. Differential Geometry, 31 (1990), pp. 215-248.
\bibitem{yau1} 
S.T. Yau, {\it On the Curvature of Compact Hermitian Manifolds}, Inv. Math. 25
 (1974), pp. 213-239.
\end{thebibliography}
\end{document}